\documentclass[a4paper,11pt]{article}
\usepackage{graphicx}
\usepackage{amsfonts}
\usepackage{amsmath}
\usepackage{amssymb}
\usepackage{indentfirst}
\usepackage{titlesec}
\usepackage{caption2}
\usepackage{color}

\def\barr{\begin{array}}
\def\earr{\end{array}}
\def\bali{\begin{aligned}}
\def\eali{\end{aligned}}
\def\bearr{\begin{eqnarray}}
\def\eearr{\end{eqnarray}}

\providecommand{\play}{\displaystyle}

\providecommand{\li}{\limits}

\providecommand{\pt}{\partial}

\providecommand{\ra}{\rightarrow}

\providecommand{\da}{\downarrow}

\providecommand{\Prob}{\mathbf P}

\providecommand{\E}{\mathbf E}

\providecommand{\al}{\alpha}

\providecommand{\bt}{\beta}

\providecommand{\Gm}{\Gamma}

\providecommand{\dt}{\delta}

\providecommand{\Dt}{\Delta}

\providecommand{\ve}{\varepsilon}

\providecommand{\tht}{\theta}

\providecommand{\kp}{\kappa}

\providecommand{\lb}{\lambda}

\providecommand{\sm}{\sigma}

\providecommand{\R}{\mathbb R}

\providecommand{\cF}{\mathcal F}

\providecommand{\iY}{\mathfrak{Y}}

\providecommand{\iiY}{\mathbf{\mathfrak{Y}}}

\providecommand{\1}{\mathbf 1}

\providecommand{\contfunc}{\mathbf{C}}

\providecommand{\boldx}{\boldsymbol{x}}

\providecommand{\boldX}{\boldsymbol{X}}

\providecommand{\boldn}{\boldsymbol{n}}

\providecommand{\boldnu}{\boldsymbol{\nu}}

\providecommand{\vphi}{\varphi}

\providecommand{\Vbar}{\overline{V}}

\providecommand{\Xhat}{\widehat{X}}

\providecommand{\Zhat}{\widehat{Z}}

\begin{document}

\title{On diffusion in narrow random channels}
\author{Mark Freidlin\thanks{Department of Mathematics,
University of Maryland at College Park, mif@math.umd.edu.} \ , \
Wenqing Hu\thanks{Department of Mathematics, University of Maryland
at College Park, huwenqing@math.umd.edu.}}

\date{}

\maketitle

\begin{abstract}

We consider in this paper a solvable model for the motion of
molecular motors. Based on the averaging principle, we reduce the
problem to a diffusion process on a graph. We then calculate the
effective speed of transportation of these motors.

\end{abstract}

\textit{Keywords:} Brownian motors/ratchets, averaging principle,
diffusion processes on graphs, random environment.

\textit{2010 Mathematics Subject Classification Numbers:} 60H30,
60J60, 92B05, 60K37.

\section{Introduction}

One of the possible ways to model Brownian motors/ratchets is to
describe them as particles (which model the protein molecules)
traveling along a designated track (see \cite{[RMP Molecular
Motors]}). At a microscopic scale such a motion is conveniently
described as a diffusion process with a deterministic drift. On the
other hand, the designated track along which the molecule is
traveling can be viewed as a tubular domain of some random shape. In
particular, such a domain can have many random "wings" added to it.
(See Fig.1. The shaded areas represent the "wings".)
 In this paper we are going to introduce a mathematically solvable model
of the Brownian motor and discuss some interesting relevant
questions around this problem. Our model is based on ideas similar
to that of \cite{[FW fish]} and \cite[Chapter 7]{[F green book]}.

The model is as follows. Let $h_0^{\pm}(x)$ be a pair of piecewise
smooth functions with $h_0^+(x)-h_0^-(x)=l_0(x)>0$. Let
$D_0=\{(x,z): x\in \R , h_0^-(x)\leq z \leq h_0^+(x)\}$ be a tubular
$2$-d domain of infinite length, i.e. it goes along the whole
$x$-axis. At the discontinuities of $h_0^{\pm}(x)$, we connect the
pieces of the boundary via straight vertical lines. The domain $D_0$
models the "main" channel in which the motor is traveling. Let a
sequence of "wings" $D_j$ ($j\geq 1$) be attached to $D_0$. These
wings are attached to $D_0$ at the discontinuities of the functions
$h_0^{\pm}(x)$.

Consider the union $D=D_0\bigcup\left(\bigcup\li_{j=1}^\infty
D_j\right)$. An example of such a domain $D$ is shown in Fig.1, in
which one can see four "wings" $D_1,D_2,D_3,D_4$. We assume that,
after adding the "wings", for the domain $D$, the boundary $\pt D$
has two smooth pieces: the upper boundary and the lower boundary.
Let $\boldn(x,z)=(n_1(x,z), n_2(x,z))$ be the inward unit normal
vector to $\pt D$. We make some assumptions on the domain $D$.

\textbf{Assumption 1.} The set of points $x\in \R$ for which there
are points $(x,z)\in \pt D$ at which the unit normal vector
$\boldn(x,z)$ is parallel to the $x$-axis: $n_2(x,z)=0$ has no limit
points in $\R$. Each such point $x$ corresponds to only one point
$(x,z)\in \pt D$ for which $n_2(x,z)=0$.

\textbf{Assumption 2.} For every $x$ the cross-section of the region
$D$ at level $x$, i.e., the set of all points belonging to $D$ with
the first coordinate equal to $x$, consists of either one or two
intervals that are its connected components. That is to say, in the
case of one interval this interval corresponds to the "main channel"
$D_0$; and in the case of two intervals one of them corresponds to
the "main channel" $D_0$ and the other one corresponds to the wing.
The wing will not have additional branching structure. Also, for
some $0<l_0<\bar{l}_0<\infty$ we have $\l_0\leq
h_0^+(x)-h_0^-(x)=l_0(x)\leq \bar{l}_0$.

Let us take into account randomness of the domain $D$. Keeping the
above assumptions in mind, we can assume that the functions
$h_0^{\pm}(x)$ and the shape of the wings $D_k$ ($k=1,2,...$) are
all random. Thus we can view the shape of $D$ as random. We
introduce a filtration $\cF_s^t$, $-\infty\leq s<t\leq \infty$ as
the smallest $\sm$-algebra corresponding to the shape of
$D\cap\{(x,z): x\in [s, t]\}$. We introduce stationarity and mixing
assumptions. Let us consider some $A\in \cF_s^t$, $-\infty \leq
s<t\leq \infty$. The set $A$ consists of some shapes of the domain
$D\cap \{(x,z): x\in [s,t]\}$. Let $\tht_r$ ($r\in \R$) be the
operator corresponding to the shift along $x$-direction:
$\tht_r(A)\in \cF_{s+r}^{t+r}$ consists of the same shapes as those
in $A$ but correspond to the domain $D\cap \{(x,z): x\in
[s+r,t+r]\}$.

\textbf{Assumption 3.} (stationarity) We have $\Prob(A)=\Prob(\tht_r(A))$.

\textbf{Assumption 4.} (mixing) For any $A\in \cF_s^t$ and any $B\in
\cF_{s+r}^{t+r}$ we have
$$\lim\li_{r\ra \pm\infty}|\Prob(A\cap B)-\Prob(A)\Prob(B)|= 0$$
exponentially fast.

For instance, we can assume that there exists some $M>0$ such that
$\Prob(A\cap B)=\Prob(A)\Prob(B)$ for $|r|\geq M$.

Here and below the symbols $\Prob$ and $\E$ etc. refer to
probabilities and expectations etc. with respect to the filtration
$\{\cF_s^t\}_{-\infty\leq s <t\leq \infty}$.

In many problems it is natural to assume that the domain $D$ is a
thin and long channel. This leads to the formulation of the problem
as follows. Let $D^\ve=\{(x,\ve z): (x,z)\in D\}$. The parameter
$\ve>0$ is small. Consider the diffusion process
$\widehat{\boldX}_t^\ve=(\Xhat_t^\ve, \Zhat_t^\ve)$ in the domain
$D^\ve$, which is described by the following system of stochastic
differential equations:

$$\left\{\begin{array}{l}
d\Xhat_t^\ve=dW_t^1+V(\Xhat_t^\ve,
\Zhat_t^\ve/\ve)dt+\nu_1(\Xhat_t^\ve,
\Zhat_t^\ve)d\widehat{\ell}_t^\ve \ ,
\\
d\Zhat_t^\ve=dW_t^2+\nu_2(\Xhat_t^\ve,
\Zhat_t^\ve)d\widehat{\ell}_t^\ve \ .
\end{array} \right.  \eqno(1)$$

Here the scalar field $V(x,z)>0 , (x,z)\in D$ characterizes the
speed of the transportation in the $x$-direction. The vector field
$\boldnu=(\nu_1,\nu_2)$ on $\pt D^\ve$ is defined as the inward unit
normal vector at the corresponding point on $\pt D$: $\boldnu(x,\ve
z)=\boldn(x,z)$ when $(x,z)\in \pt D$. The process $(W_t^1, W_t^2)$
is a standard $2$-dimensional Wiener process. We make another
assumption here.

\textbf{Assumption 5.} The process $(W_t^1, W_t^2)$ is independent
of the filtration $\{\cF_s^t\}_{-\infty\leq s<t\leq \infty}$
corresponding to the shape of $D$.

In other words our process $\widehat{\boldX}_t^\ve$ is moving in an
independent random environment characterized by random shape of the
domain $D$.

The process $\widehat{\ell}_t^\ve$ is the local time of the process
$\widehat{\boldX}_t^\ve$ at $\pt D^\ve$.

\begin{figure}
\centering
\includegraphics[height=7cm, width=10cm , bb=24 20 403 289]{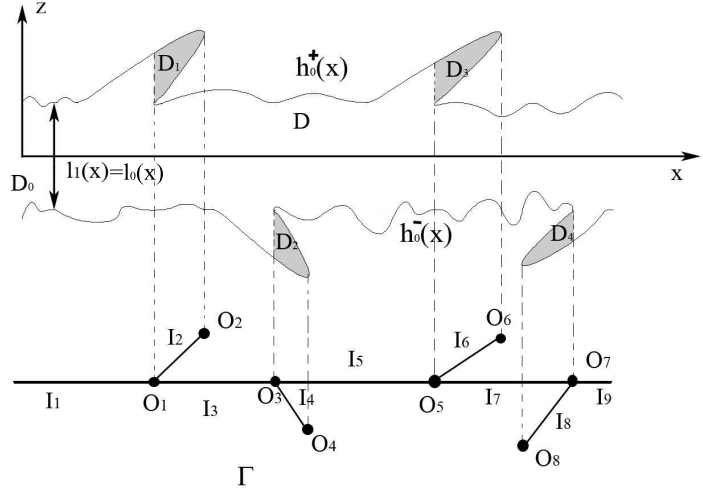}
\caption{A model of the molecular motor.}
\end{figure}

Making a change of variable $\Zhat_t^\ve\ra \Zhat_t^\ve/\ve=Z_t^\ve$
in the equation (1), we can equivalently consider the diffusion
process $\boldX_t^\ve=(X_t^\ve,Z_t^\ve)$ in the original domain $D$
as follows:

$$\left\{\begin{array}{l}
dX_t^\ve=dW_t^1+V(X_t^\ve, Z_t^\ve)dt+\nu_1^\ve(X_t^\ve,
Z_t^\ve)d\ell_t^\ve \ ,
\\
dZ_t^\ve=\dfrac{1}{\ve}dW_t^2+\nu_2^\ve(X_t^\ve, Z_t^\ve)d\ell_t^\ve
\ ,
\end{array} \right. \eqno(2)$$

Here $\boldnu^\ve$ is the co-normal vector field corresponding to
the operator $\dfrac{1}{2\ve^2}\dfrac{\pt^2}{\pt
z^2}+\dfrac{1}{2}\dfrac{\pt^2}{\pt x^2}$:
$\boldnu^\ve(x,z)=(\nu^\ve_1(x,z), \nu^\ve_2(x,z))\equiv (\ve
n_1(x,z), n_2(x,z))$. The process $\ell_t^\ve$ is the local time of
the process $\boldX_t^\ve$ at $\pt D$.

Here and below we use the symbols $\Prob^W$ and $\E^W$ etc.
(sometimes with a subscript to denote the starting point of the
process) to refer to probabilities and expectations etc. with
respect to the filtration generated by the 2-d Wiener process
$(W_t^1,W_t^2)$.

The process $\boldX_t^\ve$ has the "fast" and the "slow" components.
The "fast" component is the process $Z_t^\ve$ and the "slow"
component is the process $X_t^\ve$. According to the averaging
principle we can expect a mixing in the "fast" component before the
"slow" component $X_t^\ve$ changes significantly. We shall describe
the limiting slow motion.

In the next section we will characterize the limiting slow motion,
which is a diffusion process on a graph. This graph corresponds to
the domain $D$. A sketch of the proof of this result is in Section
3.

An interesting question arising in the applications is to calculate
the effective speed of the particles. In mathematical language this
problem can be formulated as follows. Let $\sm^\ve((-\infty,a])$ be
the first time that the process $\boldX_t^\ve$, starting from a
point $\boldx_0=(x_0,z_0)\in D$, hits $D\cap \{x=a>0\}$. The limit
$$\lim\li_{a\ra \infty}\lim\li_{\ve\da
0}\dfrac{\sm^\ve((-\infty,a])}{a}$$ exists in
$\Prob\times\Prob^W_{(x_0,z_0)}$-probability and can be viewed as
the inverse of the average effective speed of transportation of the
particle inside $D$. Using the results in Sections 2 and 3 we can
calculate this limit. This is done in Section 4. (In particular, see
Theorem 7.)

In the last Section 5 we mention briefly problems for
multidimensional channels, for random channels changing in time, and
some other generalizations.

\section{The limiting process}

Let us, for the present and for the next section, work with a fixed
shape of $D$. In the language of random motions in random
environment the convergence results that we are going to state are
in the so called "quenched" setting. We will allow this shape to be
random in Section 4. We shall find the limiting slow motion of the
diffusion process $\boldX_t^\ve$ inside $D$.

First of all we need to construct from the domain $D$ a graph $\Gm$
(see Fig.1). For $x_0\in \R$ let $C(x_0)=\{(x,z)\in D: x=x_0\}$ be
the cross-section of the domain $D$ with the line $\{x=x_0\}$. The
set $C(x_0)$ may have several connected components. We identify all
points in each connected component and the set thus obtained,
equipped with the natural topology, is homeomorphic to a graph
$\Gm$. We label the edges of this graph $\Gm$ by $I_1,...,I_k,...$
(there might be infinitely many such edges).

We see that the structure of the graph $\Gm$ consists of many edges
(such as $I_1, I_3, I_5, I_7, I_9$,... in Fig.1) that form a long
line corresponding to the domain $D_0$ and many other short edges
(such as $I_2, I_4, I_6, I_8$,... in Fig.1) attached to the long
line in a random way.

A point $y\in \Gm$ can be characterized by two coordinates: the
horizontal coordinate $x$, and the discrete coordinate $k$ being the
number of the edge $I_k$ in the graph $\Gm$ to which the point $y$
belongs. Let the identification mapping be $\iY: D \ra \Gm$. We note
that the second coordinate is not chosen in a unique way: for $y$
being an interior vertex $O_i$ of the graph $\Gm$ we can take $k$ to
be the number of any of the several edges meeting at the vertex
$O_i$.

The distance $\rho(y_1,y_2)$ between two points $y_1=(x_1, k)$ and
$y_2=(x_2, k)$ belonging to the same edge of the graph $\Gm$ is
defined as $\rho(y_1,y_2)=|x_1-x_2|$; for $y_1, y_2\in \Gm$
belonging to different edges of the graph it is defined as the
geodesic distance
$\rho(y_1,y_2)=\min(\rho(y_1,O_{j_1})+\rho(O_{j_1},
O_{j_2})+...+\rho(O_{j_l}, y_2))$, where the minimum is taken over
all chains $y_1\leftrightarrow O_{j_1}\leftrightarrow
O_{j_2}\leftrightarrow...\leftrightarrow O_{j_l}\leftrightarrow y_2$
of vertices $O_{j_i}$ connecting the points $y_1$ and $y_2$.

For an edge $I_k=\{(x,k): A_k\leq x\leq B_k\}$ we consider the
"tube" $U_k=\iiY^{-1}(I_k)\cap \{A_k\leq x\leq B_k\}$ in $D$. The
"tube" $U_k$ can be characterized by the interval $x\in [A_k,B_k]$
and the "height functions" $h_k^{\pm}(x)$: $U_k=\{(x,z): A_k\leq
x\leq B_k , h_k^-(x)\leq z \leq h_k^+(x)\}$. For $x\in [A_k, B_k]$,
we denote the set $C_k(x)$ to be the connected component of $C(x)$
that corresponds to the "tube" $U_k$: $C_k(x)=\{x\}\times [h_k^-,
h_k^{+}]$. Let $l_k(x)=h_k^+(x)-h_k^-(x)\geq 0$ for all $x\in \R$.
We notice that each $h_k^{\pm}(x) , l_k(x)$, etc. is smooth.

The vertices $O_j$ correspond to the connected components containing
points $(x,z)\in \pt D$ with $n_2(x,z)=0$. There are two types of
vertices: the interior vertices (in Fig.1 they are $O_1, O_3, O_5,
O_7$) are the intersection of three edges; the exterior vertices (in
Fig.1 they are $O_2, O_4, O_6, O_8$) are the endpoints of only one
edge.

Using the ideas in \cite{[FW fish]} with a little modification we
can establish the weak convergence of the process
$Y_t^\ve=\iiY(\boldX_t^\ve)$ (which is not Markov in general) as
$\ve \da 0$ in the space $\contfunc_{[0,T]}(\Gm)$ to a certain
Markov process $Y_t$ on $\Gm$. A sketch of the proof of this fact is
in the next section.

The process $Y_t$ is a diffusion process on $\Gm$ with a generator
$A$ and the domain of definition $D(A)$. We are going now to define
the operator $A$ and its domain of definition $D(A)$.

For each edge $I_k$ we define an operator $\overline{L}_k$:

$$\overline{L}_k u(x)=\dfrac{1}{2l_k(x)}\dfrac{d}{dx}
\left(l_k(x)\dfrac{du}{dx}\right)+\overline{V}_k(x)\dfrac{du}{dx}  \
, \ A_k \leq x \leq B_k \ . $$

Here
$$\overline{V}_k(x)=\dfrac{1}{l_k(x)}\int_{h_k^-(x)}^{h_k^+(x)}V(x,z)dz$$
is the average of the velocity field $V(x,z)$ on the connected
component $C_k(x)$, with respect to Lebesgue measure in
$z$-direction. At places where $l_k=0$, the above expression for
$\overline{V}_k(x)$ is understood as a limit as $l_k\ra 0$:
$$\overline{V}_k(x)=\lim\li_{y\ra x}\dfrac{1}{l_k(y)}\int_{h_k^-(y)}^{h_k^+(y)}V(y,z)dz \ .$$

For simplicity of presentation we will assume throughout this paper
the following.

\textbf{Assumption 6.} The function $\Vbar_k(x)=\bt>0$ is a
constant.

The case of non-constant $\Vbar_k(x)$ can be treated in a similar
way. The only difference is that the calculations are a little bit
more bulky. To be more precise, in the ordinary differential
equations we are going to solve in the proof of Theorem 2 and Lemma
1 the constant $\bt$ will be replaced by $\overline{V}_k(x)$, and
these equations can be solved correspondingly.

We also let $$\overline{L_k^0}
u(x)=\play{\dfrac{1}{2l_k(x)}\dfrac{d}{dx}\left(l_k(x)\dfrac{du}{dx}\right)}
\ .$$

The operator $\overline{L_k^0}$ can be represented as a generalized
second order differential operator (see \cite{[Feller]})

$$\overline{L_k^0} u(x)=D_{r_k}D_{q_k} f(x)  \ ,$$
where, for an increasing function $h$, the derivative $D_h$ is
defined by $D_h g(x)=\lim\li_{\dt \da
0}\dfrac{g(x+\dt)-g(x)}{h(x+\dt)-h(x)}$, and

$$q_k(x)=\int \dfrac{dx}{l_k(x)} \ , \ r_k(x)=2\int l_k(x)dx \ .$$

The operator $A$ is acting on functions $f$ on the graph $\Gm$: for
$y=(x,k)$ being an interior point of the edge $I_k$ we take
$Af(y)=\overline{L}_k f(x,k)$.

The domain of definition $D(A)$ of the operator $A$ consists of such
functions $f$ satisfying the following properties.

$\bullet$ The function $f$ must be a continuous function that is twice
continuously differentiable in $x$ in the interior part of every
edge $I_k$;

$\bullet$ There exist finite limits $\lim\li_{y\ra O_i}A f(y)$
(which are taken as the value of the function $Af$ at the point
$O_i$);

$\bullet$ There exist finite one-sided limits $\lim\li_{x\ra
x_i}D_{q_k}f(x,k)$ along every edge ending at $O_i=(x_i,k)$ and
 they satisfy the gluing conditions

$$\sum\li_{j=1}^{N_i}(\pm)\lim\li_{x\ra x_i}D_{q_{k_j}}f(x,k_j)=0 \ , \eqno(3)$$
where the sign "$+$" is taken if the values of $x$ for points
$(x,k_j)\in I_{k_j}$ are $\geq x_i$ and "$-$" otherwise. Here
$N_i=1$ (when $O_i$ is an exterior vertex) or $3$ (when $O_i$ is an
interior vertex).

For an exterior vertex $O_i=(x_i,k)$ with only one edge $I_k$
attached to it the condition (3) is just $\lim\li_{x\ra
x_i}D_{q_k}f(x,k)=0$. Such a boundary condition can also be
expressed in terms of the usual derivatives $\dfrac{d}{dx}$ instead
of $D_{q_k}$. It is $\lim\li_{x\ra x_i}l_k(x)\dfrac{\pt f}{\pt
x}(x,k)=0$. We remark that we are in dimension 2 so that these
exterior vertices are accessible, and the boundary condition can be
understood as a kind of (not very standard) instantaneous
reflection. In dimension 3 or higher these endpoints do not need a
boundary condition, they are just inaccessible. For an interior
vertex the gluing condition (3) can be written with the derivatives
$\dfrac{d}{dx}$ instead of $D_{q_k}$. For $k$ being one of the $k_j$
we define $\al_{ik}=\lim\li_{x\ra x_i}l_k(x)$ (for each edge $I_k$
the limit is a one-sided one). Then the condition (3) can be written
as

$$\sum\li_{j=1}^3 (\pm)\al_{i,k_j}\cdot\lim\li_{x\ra x_i}\dfrac{df(x,k_j)}{dx}=0 \ . \eqno(4)$$

It can be shown as in \cite[Section 2]{[FW fish]} that the process
$Y_t$ exists as a continuous strong Markov process on $\Gm$.

We fix the shape of $D$. For every $\ve>0$, every $\boldx=(x,z)\in
D$ and every $T\in (0,\infty)$ let us consider the distribution
$\mu_{\boldx}^\ve$ of the trajectory $Y_t^\ve=\iiY(\boldX_t^\ve)$
starting from a point $\boldX_0^\ve=\boldx$ in the space
$\contfunc_{[0,T]}(\Gm)$ of continuous functions on the interval
$[0,T]$ with values in $\Gm$: the probability measure defined for
every Borel subset $B\subseteq \contfunc_{[0,T]}(\Gm)$ as
$\mu_{\boldx}^\ve(B)=\Prob_{\boldX_0^\ve=\boldx}^W(Y_{\bullet}^\ve\in
B)$. Similarly, for every $y\in \Gm$ and $T>0$ let $\mu^0_y$ be the
distribution of the process $Y_t$ in the same space:
$\mu^0_y(B)=\Prob^W_y(Y_{\bullet}\in B)$. The following theorem is
our main tool for the analysis.

\textbf{Theorem 1.} \textit{For every $\boldx\in D$ and every $T>0$
the distribution $\mu_{\boldx}^\ve$ converges weakly to to
$\mu_{\iiY(\boldx)}^0$ as $\ve\da 0$.}

In other words we have $$\E_{\boldX_0^\ve=\boldx}^W
F(Y^\ve_\bullet)\ra \E_{\iiY(\boldx)}^W F(Y_{\bullet})$$ for every
bounded continuous functional $F$ on the space
$\contfunc_{[0,T]}(\Gm)$.

\section{Sketch of the proof of the convergence}

We shall now briefly give a sketch of the proof of Theorem 1
announced in the previous section.

The averaging within each edge $I_k$ is a routine adaptation of the
arguments of \cite[Section 3]{[FW fish]}. Within one edge $I_k$, the
motion of the component $X_t^\ve$ is given by the integral form of
the stochastic differential equation
$$X_t^\ve=W_t^1+\int_0^tV(X_s^\ve, Z_s^\ve)ds+\int_0^t\nu_1^\ve(X_s^\ve,Z_s^\ve)d\ell_s^\ve \ ,$$
and the one for the limiting motion $X_t$ looks like
$$X_t=W_t^1+\int_0^t \overline{V}_k(X_s)ds+\dfrac{1}{2}\int_0^t \dfrac{l_k'(X_s)}{l_k(X_s)}ds \ .$$

From the above two formulas we see that in order to prove the
convergence of $X_t^\ve$ to $X_t$ as $\ve \da 0$ in the interval
$0\leq t \leq T<\infty$ we just need the estimates of
$$(I)=\left|\int_0^t \nu_1^\ve(X_s^\ve,
Z_s^\ve)d\ell_s^\ve-\dfrac{1}{2}\int_0^t
\dfrac{l_k'(X_s^\ve)}{l_k(X_s^\ve)}ds\right|$$ and
$$(II)=\left|\int_0^t V(X_s^\ve,Z_s^\ve)ds-\int_0^t
\Vbar_k(X_s^\ve)ds\right| \ .$$

The estimate of $(I)$ ("averaging with respect to local time") is
exactly the same as that of \cite[Section 3]{[FW fish]}. For the
estimate of $(II)$ we can introduce an auxiliary function $a_k(x,z)$
satisfying the problem

$$\left\{\begin{array}{l}
\dfrac{1}{2}\dfrac{\pt^2 a_k}{\pt z^2}(x,z)=V(x,z)-\overline{V}_k(x)
\ ,
\\
\\
\dfrac{\pt a_k}{\pt z}(x,h_k^-(x))=\dfrac{\pt a_k}{\pt
z}(x,h_k^+(x))=0 \ .
\end{array}\right.$$

The solvability of this equation is guaranteed by the fact that
$$\play{\int_{h_k^-(x)}^{h_k^+(x)}(V(x,z)-\overline{V}_k(x))dz=0}$$
(this is the key point in averaging). The solution is bounded with
bounded derivatives. Applying the generalized It\^{o}'s formula (see
\cite[Section 3, equation (3.1)]{[FW fish]}) to the function $a_k$
we see that

$$\begin{array}{l}
a_k(X_t^\ve,Z_t^\ve)-a_k(X_0^\ve,Z_0^\ve)
\\
\ \ \ \ \ \play{=\int_0^t\dfrac{\pt a_k}{\pt x}(\boldX_s^\ve)dW_s^1
+\dfrac{1}{2}\int_0^t\dfrac{\pt^2 a_k}{\pt x^2}(\boldX_s^\ve)ds +
\int_0^t \dfrac{\pt a_k}{\pt x}(\boldX_s^\ve)V(\boldX_s^\ve)ds+}
\\
\ \ \ \ \ \ \ \play{\ve^{-1/2}\int_0^t\dfrac{\pt a_k}{\pt
z}(\boldX_s^\ve)dW_s^2+\ve^{-1}\int_0^t\dfrac{1}{2}\dfrac{\pt^2a_k}{\pt
z^2}(\boldX_s^\ve)ds \ .}
\end{array}$$

Multiplying both sides by $\ve$ and taking into account the problem
that $a_k$ satisfies it is immediate to get an estimate of $(II)$.
These justify the averaging within one edge $I_k$.

The gluing conditions can be obtained using the results of \cite{[FW
fish]} and the Girsanov formula. To this end one can introduce an
auxiliary process $\widetilde{\boldX}_t^\ve=(\widetilde{X}_t^\ve,
\widetilde{Z}_t^\ve)$ in $D$ via the following stochastic
differential equation

$$\left\{\begin{array}{l}
d\widetilde{X}_t^\ve=dW_t^1+\nu_1^\ve(\widetilde{X}_t^\ve,
\widetilde{Z}_t^\ve)d\widetilde{\ell}_t^\ve \ ,
\\
d\widetilde{Z}_t^\ve=\dfrac{1}{\ve}dW_t^2+\nu_2^\ve(\widetilde{X}_t^\ve,
\widetilde{Z}_t^\ve)d\widetilde{\ell}_t^\ve \ .
\end{array} \right.$$

Here $\widetilde{\ell}_t^\ve$ is the local time for the process
$\widetilde{\boldX}_t^\ve$ at $\pt D$.

This is exactly \cite[formula (1.4)]{[FW fish]}. The limiting
process within an edge is governed by $\overline{L_k^0}$. By
applying Theorem 1.2 in \cite{[FW fish]} we see that the gluing
condition is just the gluing condition in (3). On the other hand,
The measure $m^\ve$ corresponding to the process $\boldX_t^\ve$ is
related to the measure $\widetilde{m}^\ve$ corresponding to the
process $\widetilde{\boldX}_t^\ve$ in $\contfunc_{[0,T]}(D)$ via the
Girsanov formula

$$\dfrac{dm^\ve}{d\widetilde{m}^\ve}(\widetilde{\boldX}_t^\ve)=\exp\left\{\int_0^T V(\widetilde{\boldX}_s^\ve)dW_s^1-
\dfrac{1}{2}\int_0^TV^2(\widetilde{\boldX}_s^\ve)ds\right\} \ .$$

From the Girsanov formula one can show that the above density is
close to $1$ as $T$ is small. On the other hand, the process
$\boldX_t^\ve$ will spend a relatively small amount of time in a
neighborhood of the cross-section that corresponds to a branching. A
standard argument (see \cite{[FWeber]} and \cite[Appendix A.2]{[FH
LL]}) guarantees that the gluing conditions remain the same for the
process $\boldX_t^\ve$ and thus we have proved Theorem 1."

\section{Analysis of the limiting process}

The next goal is to quantify the effective speed of the motion of
the particle. As we have pointed out in Section 1, this is an
interesting question coming from applications. Our calculation in
this section will always be performed by first fix a shape of $D$
and then let the shape of $D$ be random.

\subsection{The case when there is only one edge of the graph $\Gm$}

The simplest case is that when there are no "wings" and also that
the graph $\Gm$ consists of only one edge $I_1$. Let the
corresponding $l_1(x)$, $\Vbar_1(x)=\bt>0$, etc. be defined.
Recalling our assumptions in Section 1, we see that in this case the
functions $h_0^{\pm}(x)$ do not have discontinuities and
$l_0(x)=l_1(x)$ is smooth and positive, uniformly bounded from above
and below. Consider the interval $(-\infty, a]$ for some $a>0$. Let
the process $Y_t$ start from $x=0$. Let $\tau((-\infty,a])$ be the
first exit time of $Y_t$ from $(-\infty,a]$. The random variable
$l_0(x)=l_1(x)$ is distributed according to our stationary and
mixing assumptions. We have the following.

\textbf{Theorem 2.} \textit{We have} $$\lim\li_{a\ra
\infty}\dfrac{\E^W_0\tau((-\infty, a])}{a}=2\int_0^\infty
K(t)\exp(-2\bt t)dt \ ,$$ \textit{where the function
$K(t)=\E\dfrac{l_0(s)}{l_0(s+t)}$.}

\textbf{Proof.} We see that $u(x)=\E^W_x \tau((-\infty,a])$
 is the solution of the problem

$$\left\{\begin{array}{l}
\dfrac{1}{2l_1(x)}\dfrac{d}{dx}\left(l_1(x)\dfrac{du}{dx}\right)+\bt\dfrac{du}{dx}=-1
\ ,
\\
u(-\infty)=u(a)=0 \ .
\end{array}\right.$$

The above problem can be solved explicitly. We shall first expand
the equation as

$$\dfrac{d^2u}{dx^2}+\left(\dfrac{l_1'(x)}{l_1(x)}+2\bt\right)\dfrac{du}{dx}=-2 \ .$$

Now we introduce

$$\begin{array}{l}
\play{\vphi(x)}
\\
\play{\equiv\exp\left[\int_0^x\left(\dfrac{l_1'(y)}{l_1(y)}+2\bt\right)dy\right]\dfrac{du}{dx}(x)}
\\
\play{=\dfrac{l_1(x)}{l_1(0)}\exp\left(2\bt
x\right)\dfrac{du}{dx}(x) \ .}
\end{array}$$

Using the equation that $u(x)$ satisfies it is not hard to check
that

$$\vphi'(x)=-2\dfrac{l_1(x)}{l_1(0)}\exp\left(2\bt x\right) \ .$$

Integration gives

$$\vphi(x)=-2\int_0^x \dfrac{l_1(y)}{l_1(0)}\exp\left(2\bt y\right)dy+\vphi(0) \ .$$

Thus

$$\begin{array}{l}
\play{\dfrac{du}{dx}=-\dfrac{2}{l_1(x)}\int_0^x
l_1(y)\exp\left(-2\bt(x-y)\right)dy+\vphi(0)\dfrac{l_1(0)}{l_1(x)}\exp\left(-2\bt
x\right) \ .}
\end{array}$$

Taking into account that $u(a)=0$ we see that

$$\begin{array}{l}
\play{u(x)=-2\int_a^x\dfrac{dy}{l_1(y)}\left(\int_0^yl_1(t)\exp\left(-2\bt
(y-t)\right)dt\right)+\vphi(0)\int_a^x
\dfrac{l_1(0)}{l_1(y)}\exp\left(-2\bt y\right)dy \ .}
\end{array}$$

Thus

$$u(0)=2\int_0^a \dfrac{dy}{l_1(y)}\int_0^y
l_1(t)\exp\left(-2\bt (y-t)\right)dt- C(a) \int_0^a
\dfrac{1}{l_1(y)}\exp\left(-2\bt y\right)dy \ .
$$

Here the constant $C(a)$ can be determined from the fact that
$\lim\li_{b\ra\infty}u(-b)=0$. This gives

$$\begin{array}{l}
C(a)
\\
=\lim\li_{b\ra\infty}\dfrac{\play{2\int_{-b}^a\dfrac{dy}{l_1(y)}
\left(\int_0^y l_1(t)\exp\left(-2\bt (y-t)\right)dt\right)}}
{\play{\int_{-b}^a\dfrac{1}{l_1(y)}\exp\left(-2\bt y\right)dy}}
\\
=\play{\lim\li_{b \ra \infty}
\dfrac{\dfrac{2}{l_1(-b)}\play{\int_0^{-b} l_1(t)\exp\left(2\bt
(b+t)\right)dt}} {\dfrac{1}{l_1(-b)}\play{\exp\left(2 \bt
b\right)}}}
\\
\play{=-2\int_{-\infty}^0 l_1(t)\exp\left(2\bt t\right)dt} \ .
\end{array}$$

We see from above that $|C(a)|\leq C<\infty$, where the constant $C$
is independent of $a$. This, combined with the fact that
$$\play{\int_0^a \dfrac{1}{l_1(y)}\exp\left(-2\bt y\right)dy}$$ is
uniformly bounded in $a$, show that the limit is equal to

$$\begin{array}{l}
\play{\lim\li_{a\ra \infty}\dfrac{u(0)}{a}}
\\
\play{=\lim\li_{a\ra \infty}\dfrac{2}{a}\int_0^a
\dfrac{dy}{l_1(y)}\left(\int_0^y l_1(t)\exp(-2\bt(y-t))dt\right)}
\\
\play{=\lim\li_{a\ra
\infty}\dfrac{2}{a}\int_0^a\dfrac{dy}{l_1(y)}\left(\int_0^y
l_1(y-t)\exp(-2\bt t)dt\right)}
\\
\play{=\lim\li_{a\ra\infty}\dfrac{2}{a}\int_0^a
ds\int_0^{a-s}\dfrac{l_1(s)}{l_1(s+t)}\exp(-2\bt t)dt \ .}
\end{array}$$

Let the random variable $L(s)=\play{\int_0^\infty
\dfrac{l_1(s)}{l_1(s+t)}\exp(-2\bt t)dt}$. We have $|L(s)|\leq L$ is
uniformly bounded. We fix an arbitrary $0<\mu<1$ and we have, for
any $\kp>0$ and any $0<s\leq \mu a$, there exist
$\al_0=\al_0(a,\kp,\mu)>0$ such that for any $a>\al_0$ we have

$$\left|\int_0^{a-s}\dfrac{l_1(s)}{l_1(s+t)}\exp(-2\bt t)dt
- L(s)\right|<\kp L \ .$$

Thus

$$\left|\dfrac{2}{a}\int_0^{\mu a} ds\int_0^{a-s}\dfrac{l_1(s)}{l_1(s+t)}\exp(-2\bt t)dt
- \dfrac{2}{a}\int_0^{\mu a} L(s) ds\right|<2\mu \kp L \ .$$

On the other hand, by our mixing and stationarity assumptions we see
that there exist $\al_1=\al_1(\kp,\mu)>0$ such that for any
$a>\al_1$ we have

$$\left|\dfrac{2}{a}\int_0^{\mu a} L(s) ds-2\mu \E L(s)\right|<\kp \ .$$

Therefore when $a>\max(\al_0,\al_1)$ we have

$$\left|\dfrac{2}{a}\int_0^{\mu a} ds\int_0^{a-s}\dfrac{l_1(s)}{l_1(s+t)}\exp(-2\bt t)dt
- 2\mu \E L(s)\right|<2\mu \kp L+\kp \ .$$

On the other hand we have

$$\left|\dfrac{2}{a}\int_{\mu a}^{a} ds\int_0^{a-s}\dfrac{l_1(s)}{l_1(s+t)}\exp(-2\bt t)dt
\right|\leq 2L(1-\mu) \ .$$

Thus when $a>\max(\al_0,\al_1)$ we have

$$\left|\dfrac{2}{a}\int_0^{a} ds\int_0^{a-s}\dfrac{l_1(s)}{l_1(s+t)}\exp(-2\bt t)dt
- 2 \E L(s)\right|<2\mu \kp L+\kp+4L(1-\mu) \ .$$

Since we can take $\mu$ arbitrarily close to $1$ we see from the
above estimate that we have

$$\lim\li_{a\ra\infty}\dfrac{2}{a}\int_0^a
ds\int_0^{a-s}\dfrac{l_1(s)}{l_1(s+t)}\exp(-2\bt t)dt=2\E L(s) \ .$$
Taking into account that $l_0(x)=l_1(x)$ in this case we conclude
with the statement of this theorem. $\square$

\subsection{The motion inside $D_0$}

Now we consider in more detail the case when the graph $\Gm$ does
not have any branching but it has many edges $I_1,...,I_k,...$ that
form a straight line. Let $O_1,...,O_k,...$ etc. be the
corresponding vertices. Let the corresponding cross-section width be
$l_1(x),..., l_k(x),...$ etc.. (See Fig.2.) In this case the
functions $h_0^\pm(x)$ have jumps. We can introduce a function
$l_0(x)=l_k(x)$ for $A_k\leq x\leq B_k$. The function $l_0(x)$ has
jumps at $O_i$'s and at the jumps of $l_0(x)$ we connect the pieces
of the boundary via vertical straight lines. In this way we form the
domain $D_0$ as we introduced in Section 1. We can find a family of
smooth functions $h_0^{\pm, \dt}(x)$ such that $h_0^{\pm,\dt}(x)$
converge as $\dt \da 0$ uniformly on compact subsets of $\R$ to
$h_0^{\pm}(x)$. Thus $h_0^{+,\dt}(x)-h_0^{-,\dt}(x)=l_0^\dt(x)$
converge as $\dt \da 0$ uniformly on compact subsets of $\R$ to the
function $l_0(x)$. Also, we can choose such $h_0^{\pm,\dt}(x)$ that
the domain $D_0^\dt=\{(x,z): h_0^{+,\dt}\leq z\leq h_0^{-,\dt}\}$
satisfies our assumptions. Consider the process $\boldX_t^{\ve,\dt}$
in the domain $D_0^\dt$ defined as in (2). Let
$Y_t^{\ve,\dt}=\iiY(\boldX_t^{\ve,\dt})$. For fixed $\dt>0$, the
graph $\Gm$ corresponding to $Y_t^{\ve,\dt}$ has only one edge
$I_1$. Using the results in Section 2 we see that as $\ve\da 0$ the
processes $Y_t^{\ve, \dt}$ converge weakly to $Y_t^\dt$. The process
$Y_t^\dt$ on $\Gm$ has a generator

$$\overline{L}^\dt u(x)=\dfrac{1}{2l_0^\dt(x)}\dfrac{d}{dx}\left(l_0^\dt(x)\dfrac{du}{dx}\right)
+\overline{V}_1^\dt(x)\dfrac{du}{dx} \ .$$

Here
$$\overline{V}_1^\dt(x)=\dfrac{1}{l_0^\dt(x)}\int_{h_0^{-,\dt}(x)}^{h_0^{+,\dt}(x)}V(x,z)dz$$
and we have $\overline{V}_1^\dt(x)\ra \overline{V}_k(x)=\bt>0$ as
$\dt\ra 0$ uniformly on compact subsets of $x\in \R$.

The above operator can be written in the form of a
$D_{v^\dt}D_{u^\dt}$ operator in the sense of \cite{[Feller]}. We
can explicitly calculate the functions $u^\dt$ and $v^\dt$ as
follows:

$$u^\dt(x)=\int_0^x \dfrac{1}{l_0^\dt(y)}\exp\left(-2\int_0^y \overline{V}_1^\dt(t)dt\right)dy \ ,$$

$$v^\dt(x)=2\int_0^x l_0^\dt(y)\exp\left(2\int_0^y \overline{V}_1^\dt(t)dt\right)dy \ .$$

Consider a Markov process $Y_t^{D_0}$ on the graph $\Gm$. The
process $Y_t^{D_0}$ is governed by the generalized second order
differential operator $\overline{L}=D_vD_u$ (in the sense of
\cite{[Feller]}) where

$$u(x)=\int_0^x \dfrac{1}{l_0(y)}\exp\left(-2\int_0^y \overline{V}_k(t)dt\right)dy \ ,$$

$$v(x)=2\int_0^x l_0(y)\exp\left(2\int_0^y \overline{V}_k(t)dt\right)dy \ .$$

The domain of definition of the operator $\overline{L}$ consists of
those functions $f(x)$ that are continuous and bounded, are twice
continuously differentiable and at those points $O_i$ where $l_0(x)$
have jumps it satisfies a gluing condition
$\dfrac{1}{l_{0,-}(O_i)}f'_-(O_i)=\dfrac{1}{l_{0,+}(O_i)}f'_+(O_i)$.
Here $l_{0,\pm}(O_i)$ are the corresponding values of the one
sided-limits of $l_0(x)$ to the left and to the right of $O_i$, and
$f'_{\pm}(O_i)$ the corresponding left and right derivatives. The
function $\overline{L}f(x)$ is continuous on $\Gm$.

It follows from the classical result of \cite{[FW Dynkin Fest]} that
we have the following.

\textbf{Theorem 3.} \textit{As $\dt\da 0$ the processes $Y_t^\dt$
converge weakly to the process $Y_t^{D_0}$ on
$\contfunc_{[0,T]}(\Gm)$.}

Let $\tau^\dt((-\infty,a])$ be the first time when the process
$Y_t^\dt$ exits from $(-\infty,a]$. Let $\tau^{D_0}((-\infty,a])$ be
the first time when the process $Y_t^{D_0}$ exits from
$(-\infty,a]$. Weak convergence of processes $Y_t^\dt$ as $\dt \da
0$ to $Y_t^{D_0}$ and finiteness of $\E_0^W\tau^\dt((-\infty,a])$
and $\E_0^W\tau^{D_0}((-\infty,a])$ for fixed $\dt>0$ and $a>0$
imply that we have $\lim\li_{\dt \da
0}\E_0^W\tau^\dt((-\infty,a])=\E_0^W\tau^{D_0}((-\infty,a])$.

We recall the first differential equation we used in the proof of
Theorem 2. We plug in $l_1(x)=l_0^\dt(x)$ in that equation and we
see that the corresponding solution is just the solution we get
there with $l_1(x)$ replaced by $l_0^\dt(x)$. However
$\E_0^W\tau^\dt((-\infty,a])$ is the solution of the same problem we
used in the proof of Theorem 2 with $l_1(x)$ replaced by
$l_0^\dt(x)$ and $\bt$ replaced by $\overline{V}_1^\dt(x)$. Since
$\overline{V}_1^\dt(x)\ra \bt$ as $\dt\da 0$ we see that we have

$$\begin{array}{l}
\E_0^W\tau^{D_0}((-\infty,a])=\lim\li_{\dt \da
0}\E_0^W\tau^\dt((-\infty,a])
\\
\play{=\lim\li_{\dt \da 0}\left(2\int_0^a
\dfrac{dy}{l_0^\dt(y)}\int_0^y l_0^\dt(t)\exp\left(-2\bt
(y-t)\right)dt+ 2\int_{-\infty}^0 l_0^\dt(t)\exp\left(2\bt
t\right)dt \int_0^a \dfrac{1}{l_0^\dt(y)}\exp\left(-2\bt
y\right)dy\right)}
\\
\play{=2\int_0^a \dfrac{dy}{l_0(y)}\int_0^y l_0(t)\exp\left(-2\bt
(y-t)\right)dt+2\int_{-\infty}^0 l_0(t)\exp\left(2\bt t\right)dt
\int_0^a \dfrac{1}{l_0(y)}\exp\left(-2\bt y\right)dy} \ .
\end{array}$$

Following our stationarity and mixing assumptions, after deleting
the "wings", the remaining channel still satisfies the stationarity
and mixing assumptions. Therefore by the same calculation as in the
proof of Theorem 2, and using the above formula, we see that we have
the following.

\textbf{Theorem 4.} \textit{For the process $Y_t^{D_0}$ defined as
above we have}

$$\lim\li_{a\ra \infty}\dfrac{\E_0^W\tau^{D_0}((-\infty,
a])}{a}=2\int_0^\infty K(t)\exp(-2\bt t)dt \ ,
$$ \textit{where $K(t)=\dfrac{l_0(s)}{l_0(s+t)}$ and we allow jumps of the function $l_0(x)$.
}

On the other hand, using the same argument of \cite{[FW fish]} one
can show that the process $Y_t^{D_0}$ can be viewed as the limiting
slow motion as $\ve\da 0$ of the \textit{part} of the process
$\boldX_t^\ve$ within the domain $D_0$. To be precise, consider the
domain $D$ introduced in Section 1 and the corresponding process
$\boldX_t^\ve$ in $D$. Let $\vphi_t=\play{\int_0^t
\1(\boldX_s^\ve\in D_0)ds}$ be an additive functional. (It is called
the proper time of the domain $D_0$, see \cite{[Molchanov]}.) We
introduce the time $\bt_t$ inverse to $\vphi_t$ and continuous on
the right. Let $Y_t^{\ve, D_0}=\iiY(\boldX_{\bt_t}^\ve)$. Then we
can use the same arguments of \cite{[FW fish]} to prove the weak
convergence as $\ve\da 0$ of the processes $Y_t^{\ve, D_0}$ to the
process $Y_t^{D_0}$.

Consider the process $\boldX_t^\ve$ moving in the domain $D$ as in
Section 1. Let $l_0(x)$ be the cross-section width corresponding to
the domain $D_0$. As before we see that it could have jumps. Let
$Y_t$ be the limiting slow motion as $\ve\da 0$ of $\boldX_t^\ve$,
defined as in Section 2. Let $\tau((-\infty,a])$ be the first time
the process $Y_t$ exits from $(-\infty,a]$. Since $Y_t$ is the
limiting slow motion of $\boldX_t^\ve$ and $Y_t^{D_0}$ is the
limiting slow motion of the part of $\boldX_t^\ve$ inside $D_0$, we
see from the above discussions that we have
$$\tau^{D_0}((-\infty,a])=\int_0^{\tau((-\infty,a])}\1(\iiY^{-1}(Y_t)\subset D_0)dt \ .$$

By Theorem 4 and the above relation we see that we have the
following.

\textbf{Corollary 1.} $$\lim\li_{a\ra
\infty}\dfrac{\E_0^W\play{\int_0^{\tau((-\infty,a])}\1(\iiY^{-1}(Y_t)\subset
D_0)dt}}{a}=2\int_0^\infty K(t)\exp(-2\bt t)dt$$ \textit{where}
$K(t)=\E\dfrac{l_0(s)}{l_0(s+t)}$.

This fact will be used in the next subsection.

\begin{figure}
\centering
\includegraphics[height=5cm, width=10cm , bb=47 67 408 256]{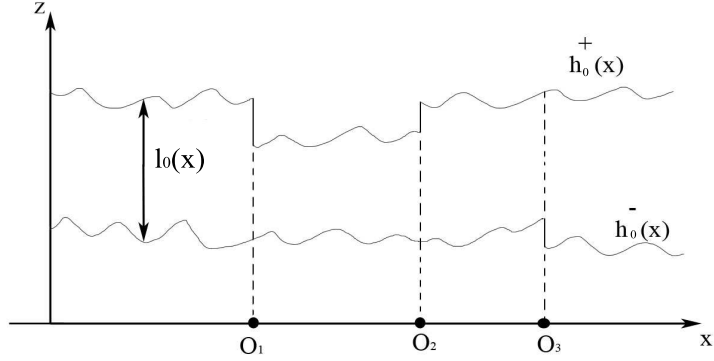}
\caption{The case when $l_0(x)$ has jumps.}
\end{figure}

\subsection{The general case}

In the general case when there is branching the domain $D$ consists
of a domain $D_0$ that has a cross-section width $l_0(x)$ which
allows occasional jumps. The "wings" $D_j$ ($j\geq 1$) are then
attached at the jumps of the domain $D_0$.

In this case our graph $\Gm$ consists of two types of edges. The
first type of edges correspond to the domain $D_0$ as we discussed
in the previous section. The second type of edges correspond to the
"wings" attached to $D_0$.

In order to calculate the effective speed of transportation we shall
first calculate the expected time that the process spends at one
fixed "wing". As a first step we do not consider the random shape
but perform the calculation for a fixed shape. Also, we shall first
consider the simplest case that $\Gm$ has only three edges:
$I_1=(-\infty,0]$, $I_2=[0,\infty)$ and $I_3$ is an edge with one
endpoint $O_1=0$ and another endpoint $O_2=(r,3)$. In this case
$l_1(x), l_2(x)$ and $l_3(x)$ are smooth functions.

We construct the process $Y_t$ as in Section 2 corresponding to the
above graph $\Gm$. Consider the interval $(-\infty, a]$ for some
$a>0$. Let the process $Y_t$ start from the point $O_1=0$. Let
$\tau_{I_3}((-\infty,a])$ be the time that the process $Y_t$ spends
at the edge $I_3$ before its first exit time from $(-\infty,a]$. We
have the following.

\textbf{Lemma 1.}

$$\E_0^W \tau_{I_3}((-\infty,a])=2\text{sign}(r)
\int_0^rl_3(t)\exp(2\bt t)dt \int_0^a \dfrac{1}{l_2(y)}\exp(-2\bt
y)dy \ .$$

\textbf{Proof.} Let for example $r>0$. The function
$u(x,k)=\E^W_{(x,k)} \tau_{I_3}((-\infty,a])$ is the solution of the
problem

$$\left\{\begin{array}{l}
\dfrac{1}{2l_k(x)}\dfrac{d}{dx}\left(l_k(x)\dfrac{du}{dx}\right)+\bt\dfrac{du}{dx}=0
\ , \text{ for } (x,k)\in I_k \text{ and } k=1,2 \ ,
\\
\dfrac{1}{2l_3(x)}\dfrac{d}{dx}\left(l_3(x)\dfrac{du}{dx}\right)+\bt\dfrac{du}{dx}=-1
\ , \text{ for } (x,3)\in I_3 \ ,
\\
u(-\infty,1)=u(a,2)=0 \ ,
\\
\lim\li_{y=(x,3)\ra O_2}l_3(x)\dfrac{du}{dx}((x,3))=0 \ ,
\\
l_1(0) \lim\li_{x\ra O_1}\dfrac{du}{dx}(x,1)=l_2(0) \lim\li_{x\ra
O_1}\dfrac{du}{dx}(x,2)+l_3(0)\lim\li_{x\ra O_1}\dfrac{du}{dx}(x,3)
\ .
\end{array}\right.$$

The function $u(x,k)$ is continuous at point $O_1$.

Similarly as in Section 4.1 there are solutions $u_1(x)=u(x,1),
u_2(x)=u(x,2)$ and $u_3(x)=u(x,3)$ corresponding to the edges $I_1$,
$I_2$ and $I_3$. They are defined as follows. We have

$$\begin{array}{l}
\play{u_1(x)=C_1l_1(0)\int_{-b}^x \dfrac{1}{l_1(t)}\exp(-2\bt
t)dt+u_1(-b)
 \text{ for } x\in (-\infty, 0] \ ,}
\end{array}$$

$$\begin{array}{l}
\play{u_2(x)=C_2l_2(0)\int_a^x \dfrac{1}{l_2(y)}\exp(-2\bt y)dy \
\text{ for } x\in [0,+\infty) \ ,}
\end{array}$$

$$\begin{array}{l}
\play{u_3(x)=-2\int_r^x\dfrac{dy}{l_3(y)}\left(\int_0^yl_3(t)\exp\left(2\bt(t-y)\right)dt
\right)+C_3l_3(0)\int_r^x\dfrac{1}{l_3(y)}\exp\left(-2\bt
y\right)dy+u_3(r) \ }
\\
\ \ \ \ \ \ \ \ \ \ \ \ \ \ \ \ \ \ \ \ \ \ \ \ \ \ \ \ \ \ \ \ \ \
\ \ \ \ \ \ \ \ \ \ \ \ \ \ \ \ \ \ \ \ \ \ \ \ \ \ \ \ \ \ \ \ \ \
\ \ \ \ \ \ \ \ \ \ \ \ \ \ \ \ \ \ \text{ for } x \in [0,r] \ .
\end{array}$$

We also know that $\lim\li_{b \ra\infty}u_1(-b)=0$. There are $4$
undetermined constants: $C_1, C_2, C_3$ and $u_3(r)$.
 We can uniquely
determine them by the $4$ relations
$$u_1(0)=u_2(0)=u_3(0) \ , $$ $$\lim\li_{x\ra r}l_3(x)\dfrac{\pt u_3}{\pt x}(x)=0 \ ,$$ and
$$-l_1(0)\dfrac{\pt u_1}{\pt x}(0)+l_2(0)\dfrac{\pt u_2}{\pt
x}(0)+l_3(0)\dfrac{\pt u_3}{\pt x}(0)=0 \ .$$

We shall solve the above system. The relations we need are as
follows

$$
\begin{array}{l}
\play{C_1l_1(0)\int_{-b}^0 \dfrac{1}{l_1(y)}\exp(-2\bt y)dy+u_1(-b)}
\\
\play{=-C_2l_2(0)\int_0^a \dfrac{1}{l_2(y)}\exp(-2\bt y)dy}
\\
\play{=-2\int_r^0 \dfrac{dy}{l_3(y)}\left(\int_0^y
l_3(t)\exp(2\bt(t-y))dt\right)+C_3l_3(0)\int_r^0
\dfrac{1}{l_3(y)}\exp(-2\bt y)dy+u_3(r) \ ,}
\end{array}
$$

$$-2\int_0^r l_3(t)\exp(2\bt(t-r))dt+C_3 l_3(0)\exp(-2\bt r)=0 \ ,$$

$$-C_1l_1(0)+C_2l_2(0)+C_3l_3(0)=0 \ .$$

Thus $$C_3l_3(0)=2\int_0^r l_3(t)\exp(2\bt t)dt \ .$$

On the other hand, we have $$C_1l_1(0)=C_2l_2(0)+C_3l_3(0) \ .$$

So we get

$$
(C_2l_2(0)+C_3l_3(0))\int_{-b}^0 \dfrac{1}{l_1(y)}\exp(-2\bt
y)dy+u_1(-b) =-C_2l_2(0)\int_0^a \dfrac{1}{l_2(y)}\exp(-2\bt y)dy \
.
$$

This gives us $\lim\li_{b \ra \infty}(C_2l_2(0)+C_3l_3(0))=0$. Thus

$$u_2(0)=2\int_0^rl_3(t)\exp(2\bt t)dt \int_0^a \dfrac{1}{l_2(y)}\exp(-2\bt y)dy \
.$$

In the case when $r<0$ the gluing condition becomes $$l_1(0)
\lim\li_{x\ra O_1}\dfrac{du}{dx}(x,1)+l_3(0) \lim\li_{x\ra
O_1}\dfrac{du}{dx}(x,3)=l_2(0)\lim\li_{x\ra O_1}\dfrac{du}{dx}(x,2)
\ .$$ Thus $C_1l_1(0)+C_3l_3(0)=C_2l_2(0)$. A similar calculation
shows that
$$\E_0^W\tau_{I_3}((-\infty,a])=2\int_r^0l_3(t)\exp(2\bt
t)dt \int_0^a \dfrac{1}{l_2(y)}\exp(-2\bt y)dy \ .$$ Thus in general
we have $$\E_0^W \tau_{I_3}((-\infty,a])=2\text{sign}(r)
\int_0^rl_3(t)\exp(2\bt t)dt \int_0^a \dfrac{1}{l_2(y)}\exp(-2\bt
y)dy \ .$$ $\square$

The above lemma can help us to deal with a more general case. Let
$q\neq 0$ and $q\in (-\infty,a]$. We assume that the graph $\Gm$
still consists of $3$ edges $I_1=(-\infty,q]$, $I_2=[q,\infty)$ and
$I_3$. In this case the edge $I_3$ has one endpoint $O_1$ lying on
$I_1\cup I_2$ but with $x$-coordinate $q$, and another endpoint
$O_2=(q+r,3)$. Let the process $Y_t$ again start from the point $0$.
(In this case the point $0$ is lying on either $I_1$ or $I_2$ but
may not be on their intersection.) Let $\tau_{I_3}((-\infty,a])$ be
the time that the process $Y_t$ spends at $I_3$ before it exits from
$(-\infty,a]$. We have the following.

\textbf{Corollary 2.} \textit{If $q>0$ then }$$\E_0^W
\tau_{I_3}((-\infty,a])=2\text{sign}(r) \int_q^{q+r}l_3(t)\exp(2\bt
(t-q))dt \int_q^a \dfrac{1}{l_2(y)}\exp(-2\bt (y-q))dy \ .$$
\textit{If $q<0$ then }$$\E_0^W
\tau_{I_3}((-\infty,a])=2\text{sign}(r) \int_q^{q+r}l_3(t)\exp(2\bt
(t-q))dt \int_0^a \dfrac{1}{l_2(y)}\exp(-2\bt (y-q))dy \ .$$

\textbf{Proof.} The above results can be easily seen from the strong
Markov property of the process $Y_t$ and Lemma 1. To be more
precise, if $q>0$, then by strong Markov property $\E_0^W
\tau_{I_3}((-\infty,a])$ is just equal to $\E_q^W
\tau_{I_3}((-\infty,a])$, which can be calculated using Lemma 1 and
a shift. If $q<0$ we need a little bit more argument. In this case
we can consider the process $Y_t$ starting from $x=q$ and its first
time of exiting from $(-\infty,0]$. Let $\tau_{I_3}((-\infty, 0])$
be the time that $Y_t$ spends at $I_3$ before it exits from
$(-\infty,0]$. Then $\E_q^W\tau_{I_3}((-\infty, 0])$ can be
calculated using Lemma 1 and a shift. On the other hand,
$\E_q^W\tau_{I_3}((-\infty,a])$ can also be calculated using Lemma 1
and a shift. From the strong Markov property of $Y_t$ we see that
$\E_0^W\tau_{I_3}((-\infty,a])=\E_q^W\tau_{I_3}((-\infty,a])-\E^W_q\tau_{I_3}((-\infty,0])$,
which gives the formula we need. $\square$

Following a similar approximation argument as we did in Section 4.2,
we can consider the case when the graph $\Gm$ consists of many edges
$I_1,I_2,...$ etc. that are of the first type. They correspond to
the domain $D_0$. We allow jumps of the function $l_0(x)$. Then we
attach only one "wing" $D_1$ to $D_0$ and the domain $D_1$
corresponds to an edge $J$ in the graph $\Gm$. Let
$l_{\text{wing}}(x)$, $r$ etc. be the quantities corresponding to
the "wing" $D_1$. Let the edge $J$ have two endpoints: $O^J_1$ lying
on $\cup_{k=1}^\infty I_k$ and with $x$-coordinate $q$; $O^J_2$
lying on the other endpoint of $J$ with $x$-coordinate $q+r$. After
an approximation and averaging we get a process $Y_t^{D_0\cup D_1}$
on $\Gm$. Let positive $a>q$. Let $\tau_{J}((-\infty,a])$ be the
time that the process $Y_t^{D_0\cup D_1}$ spends in the edge $J$
before it exits from $(-\infty,a]$. Lemma 1 and Corollary 2 and the
same approximation argument as in Section 4.2 give us the following
corollary.

\textbf{Corollary 3.} \textit{If $q>0$ then }$$\E_0^W
\tau_{J}((-\infty,a])=2\text{sign}(r)
\int_q^{q+r}l_{\text{wing}}(t)\exp(2\bt (t-q))dt \int_q^a
\dfrac{1}{l_0(y)}\exp(-2\bt (y-q))dy \ .$$\textit{ If $q<0$ then
}$$\E_0^W \tau_{J}((-\infty,a])=2\text{sign}(r)
\int_q^{q+r}l_{\text{wing}}(t)\exp(2\bt (t-q))dt \int_0^a
\dfrac{1}{l_0(y)}\exp(-2\bt (y-q))dy \ .$$ \textit{If $q=0$ then
}$$\E_0^W \tau_{J}((-\infty,a])=2\text{sign}(r)
\int_0^rl_{\text{wing}}(t)\exp(2\bt t)dt \int_0^a
\dfrac{1}{l_0(y)}\exp(-2\bt y)dy  \ .$$

Finally we come to the original problem in Section 1. We consider
the case when there are many random "wings" attached to $D_0$. Let
the process corresponding to the domain $D$ be $\boldX_t^\ve$. Let
the limiting slow motion be $Y_t$. We introduce the corresponding
quantities $l_{\text{wing}}(x)$, $r$ etc. Let
$l_0(x)=h_0^+(x)-h_0^-(x)$ be the cross section width that
corresponds to the domain $D_0$. We notice that it has occasional
jumps. Let there be $n$ wings in the interval $x\in [0,1]$. We
assume that we have the following.

\textbf{Assumption 7.} The random variable $n\geq 1$ is a bounded
random variable: $n\leq n_0<\infty$.

For $a>0$ we define the random variable

$$
M(q,r,a)=\left\{\begin{array}{l} \play{2\text{sign}(r)
\int_q^{q+r}l_{\text{wing}}(t)\exp(2\bt (t-q))dt \int_q^a
\dfrac{1}{l_0(y)}\exp(-2\bt (y-q))dy \ \text{ when } q>0 \ ;}
\\
\play{2\text{sign}(r) \int_q^{q+r}l_{\text{wing}}(t)\exp(2\bt
(t-q))dt \int_0^a \dfrac{1}{l_0(y)}\exp(-2\bt (y-q))dy \ \text{ when
} q<0 \ ;}
\\
\play{2\text{sign}(r) \int_0^rl_{\text{wing}}(t)\exp(2\bt t)dt
\int_0^a \dfrac{1}{l_0(y)}\exp(-2\bt y)dy  \text{ when } q=0 } \ .
\end{array}\right.$$

We assume that, for some constant $0<A_1<\infty$, we have the
following.

\textbf{Assumption 8.} $0\leq l_{\text{wing}}(x)\leq A_1$ and
$|r|\leq A_1$.

Thus we see that for some constant $A>0$ we have
$$\left|\int_0^r l_{\text{wing}}(t)\exp(2\bt
t)dt\right|\leq A<\infty \ .$$

We also see that $$|M(q,r,a)|\leq 2A \int_0^\infty\dfrac{1}{l_0}
\exp(-2\bt y)dy\leq M<\infty$$ for some constant $M>0$.

All these random quantities are distributed according to our
stationarity and mixing assumptions. Let $\tau((-\infty, a])$ be the
first time the process $Y_t$ exits from $(-\infty,a]$.

By our assumption on stationarity and mixing, using Corollary 2, we
see that we have the following.

\textbf{Lemma 2.}\textit{ We have}
$$\lim\li_{a\ra
\infty}\dfrac{\E_0^W\play{\int_0^{\tau(-\infty,a]}\1(\iiY^{-1}(Y_t)\not
\subset D_0)dt}}{a}=2\E n\E
\text{sign}(r)\int_0^rl_{\text{wing}}(t)\exp(2\bt t)dt \int_0^\infty
\dfrac{1}{l_0(y)}\exp(-2\bt y)dy$$ \textit{in probability.}

\textbf{Proof.} Let the "wings" located to the left of the point $0$
have $x$-coordinate $0>q_{-1}>q_{-2}>...$. Let the "wings" located
to the right of the point $0$ (including possibly the point $0$) and
not exceeding $x=a$ have $x$-coordinate $0\leq q_1<...<q_{n(a)}\leq
a$. We see that we have

$$\E_0^W \int_0^{\tau(-\infty,a]}\1(\iiY^{-1}(Y_t)\not
\subset D_0)dt=\sum\li_{k=1}^\infty
M(q_{-k},r_{-k},a)+\sum\li_{k=1}^{n(a)}M(q_k,r_k,a) \ .$$

For $q<0$ we have

$$M(q,r,a)\leq\exp(2\bt q)M(0,r,a)\leq M\exp(2\bt q) \ .$$

Thus

$$\sum\li_{k=1}^\infty M(q_{-k},r_{-k},a)\leq M\sum\li_{k=1}^\infty \exp(2\bt q_{-k}) \ .$$

Taking into account Assumption 7 we see that we have
$$\sum\li_{k=1}^\infty \exp(2\bt q_{-k})\leq n_0 \sum\li_{k=0}^\infty \exp(-2\bt k)<\infty \ .$$

Therefore

$$\lim\li_{a\ra \infty}\dfrac{\play{\E_0^W \int_0^{\tau(-\infty,a]}\1(\iiY^{-1}(Y_t)\not
\subset D_0)dt}}{a}=\lim\li_{a\ra
\infty}\dfrac{\sum\li_{k=1}^{n(a)}M(q_k,r_k,a)}{a}=\lim\li_{a\ra\infty}
\dfrac{\sum\li_{k=1}^{n(a)}M(q_k,r_k,a)}{n(a)}\dfrac{n(a)}{a} \ .$$

On the other hand, we have

$$\begin{array}{l}
\play{\sum\li_{k=1}^{n(a)}M(q_k,r_k,a)}
\\
\play{=\sum\li_{k=1}^{n(a)}2\text{sign}(r_k)
\int_{q_k}^{q_k+r_k}l_{\text{wing}}(t)\exp(2\bt
(t-q_k))dt\int_{q_k}^{a}\dfrac{1}{l_0(y)}\exp(-2\bt(y-q_k))dy}
\\
\play{=2\sum\li_{k=1}^{n(a)}W(q_k,r_k)C(q_k,a) \ .}
\end{array}$$

Here $$W(q_k,r_k)=\text{sign}(r_k)
\int_{q_k}^{q_k+r_k}l_{\text{wing}}(t)\exp(2\bt (t-q_k))dt$$ and
$$C(q_k,a)=\int_{q_k}^{a}\dfrac{1}{l_0(y)}\exp(-2\bt(y-q_k))dy \ .$$

Thus we can write
$$\sum\li_{k=1}^{n(a)}M(q_k,r_k,a)=2\sum\li_{k=1}^{n(a)}
W(q_k,r_k)C(q_k,a+q_k)-2\sum\li_{k=1}^{n(a)}W(q_k,r_k)(C(q_k,a+q_k)-C(q_k,a)) \ .$$

By the remark after Assumption 8 we see that

$$\begin{array}{l}
\play{2\sum\li_{k=1}^{n(a)}W(q_k,r_k)(C(q_k,a+q_k)-C(q_k,a))}
\\
\play{\leq 2A\sum\li_{k=1}^{n(a)}(C(q_k,a+q_k)-C(q_k,a))}
\\
\play{\leq
\dfrac{2A}{l_0}\sum\li_{k=1}^{n(a)}\int_{a}^{a+q_k}\exp(-2\bt
(y-q_k))dy }
\\
\play{=\dfrac{A}{\bt l_0}
\sum\li_{k=1}^{n(a)}[\exp(-2\bt(a-q_k))-\exp(-2\bt a)]}
\\
\play{=\dfrac{A}{\bt l_0}
\sum\li_{k=1}^{n(a)}\exp(-2\bt(a-q_k))-\dfrac{A}{\bt
l_0}n(a)\exp(-2\bt a) \ .}
\end{array}$$

Thus by our Assumption 7 again we see that
$$\lim\li_{a\ra\infty}\dfrac{2\sum\li_{k=1}^{n(a)}W(q_k,r_k)(C(q_k,a+q_k)-C(q_k,a))}{n(a)}=0 \ .$$

Therefore using the weak Law of Large Numbers for triangular arrays
and taking into account our assumptions on mixing and stationarity
we see that

$$\begin{array}{l}
\play{\lim\li_{a\ra\infty}\dfrac{\sum\li_{k=1}^{n(a)}M(q_k,r_k,a)}{n(a)}}
\\
\play{=2\lim\li_{a\ra\infty}\dfrac{\sum\li_{k=1}^{n(a)}W(q_k,r_k)C(q_k,a+q_k)}{n(a)}}
\\
\play{=2\E \text{sign}(r)\int_0^rl_{\text{wing}}(t)\exp(2\bt t)dt
\int_0^\infty \dfrac{1}{l_0(y)}\exp(-2\bt y)dy}\end{array}$$ in
probability.

To be more precise, we write

$$\begin{array}{l}
\play{\dfrac{1}{n(a)}\sum\li_{k=1}^{n(a)} M(q_k, r_k, a)}
\\
\play{=\dfrac{1}{n(a)}\sum\li_{k=1}^{n(a)}(M(q_k, r_k, a)-\E M(q_k,
r_k, a))+\E M(q_k,r_k,a) \ .}
\end{array}$$

Thus we have

$$\begin{array}{l}
\play{\lim\li_{a\ra\infty}\left(\dfrac{1}{n(a)}\sum\li_{k=1}^{n(a)}
M(q_k, r_k, a)\right)}
\\
\play{=\lim\li_{a\ra\infty}\left(\dfrac{1}{n(a)}\sum\li_{k=1}^{n(a)}(M(q_k,
r_k, a)-\E M(q_k, r_k, a))\right)+\E M(q_k,r_k,\infty) \ .}
\end{array}$$

Let $\Phi_k(a)=M(q_k,r_k,a)-\E M(q_k, r_k, a)$. We see that it
suffices to prove
$$\lim\li_{a\ra\infty}\dfrac{1}{n(a)}\sum\li_{k=1}^{n(a)}\Phi_k(a)=0$$
in probability.

Pick any $\dt>0$. By Chebyshev inequality we have the estimate

$$\begin{array}{l}
\play{\Prob\left(\left|\dfrac{1}{n(a)}\sum\li_{k=1}^{n(a)}\Phi_k(a)\right|>\dt\right)}
\\
\play{\leq \dfrac{1}{\dt^2
n^2(a)}\E\left(\sum\li_{k=1}^{n(a)}\Phi_k(a)\right)^2}
\\
\play{=\dfrac{1}{\dt^2
n^2(a)}\left(\sum\li_{k=1}^{n(a)}\E\Phi_k^2(a)+\sum\li_{k,l=1, k\neq
l}^{n(a)}\E\Phi_k(a)\Phi_l(a)\right)}
\\
\play{\leq\dfrac{1}{\dt^2 n^2(a)}\left(C_1 n(a)+C_2\sum\li_{k,l=1,
k\neq l}^{n(a)}\exp(-\lb |q_k-q_l|)\right)} \ .\end{array}$$

Here by the stationarity assumption we see that $C_1=\E\Phi_k^2(a)$
is a constant. Also, by the exponentially mixing condition we see
that
$|\E\Phi_k(a)\Phi_l(a)|=|\E\Phi_k(a)\Phi_l(a)-\E\Phi_k(a)\E\Phi_l(a)|\leq
C_2\exp(-\lb |q_k-q_l|)$, $\lb>0$. By our Assumption 7 we see that
$$\sum\li_{k,l=1, k\neq l}^{n(a)}\exp(-\lb |q_k-q_l|)\leq
\sum\li_{k,l=1, k\neq l}^{n(a)}\exp(-\dfrac{\lb}{n_0} |k-l|)<C_3
n(a)$$ for some $C_3>0$.

Thus for any $\dt>0$ we have

$$\lim\li_{a\ra\infty}\Prob\left(\left|\dfrac{1}{n(a)}\sum\li_{k=1}^{n(a)}\Phi_k(a)\right|>\dt\right)=0 \ .$$

On the other hand we have
$\play{\lim\li_{a\ra\infty}\dfrac{n(a)}{a}}=\E n$ almost surely.
Thus we can conclude with the final result. $\square$

Adding the two equations in Corollary 1 and Lemma 2 we have the
following.

\textbf{Theorem 5.} \textit{We have}

$$\begin{array}{l}
\lim\li_{a\ra \infty}\dfrac{\E_0^W\tau((-\infty,a])}{a}
\\
\\
\play{=2\int_0^\infty K(t)\exp(-2\bt t)dt +2\E n\E
\text{sign}(r)\int_0^rl_{\text{wing}}(t)\exp(2\bt t)dt \int_0^\infty
\dfrac{1}{l_0(y)}\exp(-2\bt y)dy}  \end{array}$$ \textit{in
probability. Here $K(t)=\E\dfrac{l_0(s)}{l_0(s+t)}$.}

Let $\sm^\ve((-\infty,a])$ be the first time that the process
$\boldX_t^\ve$, starting from a point $\boldx_0=(x_0,z_0)$ inside
$D$, hits the level curve $C(a)$. Since $x_0$ is finite we see that
Theorems 1 and 5 lead to the following theorem.

\textbf{Theorem 6.} \textit{We have}

$$\begin{array}{l}
\lim\li_{a\ra \infty}\lim\li_{\ve\da
0}\dfrac{\E_{\boldx_0}^W\sm^\ve((-\infty,a])}{a}
\\
\\
\play{=2\int_0^\infty K(t)\exp(-2\bt t)dt +2\E n\E
\text{sign}(r)\int_0^rl_{\text{wing}}(t)\exp(2\bt t)dt \int_0^\infty
\dfrac{1}{l_0(y)}\exp(-2\bt y)dy} \end{array}$$ \textit{in
probability. Here $K(t)=\E\dfrac{l_0(s)}{l_0(s+t)}$.}

The above theorem helps us to conclude in our original problem. We
can divide the domain $D\cap \{(x,z): x\in (-\infty,a]\}$ into
consecutive pieces alternatively of $x$-length $b$  and $\sqrt{b}$
for some $b>0$. The total time spent by the process $\boldX_t^\ve$
before it exits from $x=a$ is the sum of those times spent in
domains of $x$-length $b$ and those times spent in domains of
$x$-length $\sqrt{b}$. As we are taking $a\ra\infty$ we can also let
$b\ra\infty$ and the average time spent in domains of $x$-length
$\sqrt{b}$ will not contribute. On the other hand, since the process
$\boldX_t^\ve$ has a deterministic positive drift in the
$x$-direction, we see that as $b\ra\infty$ the motion inside
different domains of $x$-length $b$ will be asymptotically
independent. (The motion against the flow is a large deviation
effect.) Thus $\sm^\ve((-\infty,a])$ will be asymptotically
distributed as the sum of some independent random times. This leads
to the relation $\lim\li_{a\ra \infty}\lim\li_{\ve\da
0}\dfrac{\E_{\boldx_0}^W\sm^\ve((-\infty,a])}{a}=\lim\li_{a\ra\infty}\lim\li_{\ve\da
0}\dfrac{\sm^\ve((-\infty,a])}{a}$. Thus from Theorem 6 we see that
we have the following theorem.

\textbf{Theorem 7.} \textit{We have}

$$\begin{array}{l}
\lim\li_{a\ra \infty}\lim\li_{\ve\da
0}\dfrac{\sm^\ve((-\infty,a])}{a}
\\
\\
\play{=2\int_0^\infty K(t)\exp(-2\bt t)dt +2\E n\E
\text{sign}(r)\int_0^rl_{\text{wing}}(t)\exp(2\bt t)dt \int_0^\infty
\dfrac{1}{l_0(y)}\exp(-2\bt y)dy} \end{array}$$ \textit{in
probability. Here $K(t)=\E\dfrac{l_0(s)}{l_0(s+t)}$.}

\section{Remarks and generalizations}

1. The results of previous sections can be extended to the case of a
multidimensional channel $G=\{(x,z): x\in \R^1, z\in G_x\}$ where
$G_x\subset \R^{n-1}$, $x\in \R^1$ are $(n-1)$-dimensional domains
assumed to be bounded and consisting of a finite number of connected
components; each $G_x$ contains the origin. The boundary $\pt G$ of
$G\subset \R^n$ is assumed to be smooth, except, maybe, a number of
$(n-2)$-dimensional manifolds (like the discontinuity points in the
$2$-dimensional case). Let $\Gm$ be the graph homeomorphic (in the
natural topology) to the set of connected components of the sets
$\{z\in \R^{n-1}: (x,z)\in G\}$ for various $x\in \R^1$. Let the
functions $l_i(x)$ be defined as $(n-1)$-dimensional volumes of
corresponding connected components.

Consider the process $\boldX_t^\ve$ governed by the operator
$b\dfrac{\pt}{\pt x}+\dfrac{1}{2}\Dt_{x,z}$ inside the domain
$G^\ve=\{(x,z): x\in \R^1, z\in \R^{n-1}, (x,z/\ve)\in G\}$. Let
$\sm_a^\ve=\min\{t: \boldX_t^\ve\in \{(x,z)\in G\cup \pt G,
x=a\}\}$. Then, under mild additional conditions the limit
$$\lim\li_{a\ra\infty}\lim\li_{\ve\da 0}\dfrac{\sm^\ve_a}{a}$$
exists in probability $\Prob\times \Prob_{(x,z)}^W$, $(x,z)\in
G^\ve$ and is given by Theorem 7.

2. Let $\nu_t$ be a continuous time Markov chain with $2$ states $1$
and $2$. Let functions $h_i^{\pm}(x)$, $x\in \R^1$, $i=1,2$, be
piecewise smooth, and $h_i^-(x)<0<h_i^+(x)$, $x\in \R^1$. Put
$D_i^\ve=\{(x,z): x\in \R^1, -h_i^-(x)\leq z/\ve \leq h_i^+(x)\}$
and $D^\ve(t)=D_{\nu(t)}^\ve$. Define the process $(X_t^\ve,
Z_t^\ve)$ in $D^\ve(t)$ as the process governed by the operator
$b\dfrac{\pt}{\pt x}+\dfrac{1}{2}\Dt_{x,z}$ inside $D^\ve(t)$ with
the normal reflection on the boundary at the times when $\nu_t$ is
continuous (then $\nu_t$ is constant). Let $Z_t^\ve$ jump to $0$ at
times when $\nu_t$ has jumps. (Actually, we need this condition to
define the process $(X_t^\ve, Z_t^\ve)$ in a unique way; it is not
important since we are interested in the limit as $\ve\da 0$.) Then
one can prove that the slow component $X_t^\ve$ of the process
converges as $\ve\da 0$ to the process described by the equation
$\dot{X}_t=\dot{W}_t+\dfrac{1}{2}\widetilde{b}(t,X_t)$ where
$\widetilde{b}(t,x)=b+\dfrac{d}{dx}(\ln l_{\nu_t}(x))$,
$l_{\nu_t}(x)=h_{\nu_t}^+(x)-h_{\nu_t}^-(x)$. It is known that for
certain drift terms $\widetilde{b}(t,x)$, the process $X_t$
demonstrates the so called ratchet effect: If $\nu_t$ is identically
equal to $1$ or $2$, the process tends to $+\infty$ as $t\ra
\infty$, but if $\nu_t$ is the Markov chain (independent of $W_t$)
the process $X_t$ tends to $-\infty$.

One can conclude from our considerations, that the ratchet effect
can be caused by random and independent of the basic Brownian motion
changes of the geometry of the domain.

3. We assumed that the "wings" have a simple structure -- each of
them corresponds to just one edge of the graph (Fig.1). One can
consider the case of more complicated wings, like, for instance, at
the vertex $O_1$ in Fig.3. One can also include in the consideration
the case when the "obstacles" in the channel are such that the
corresponding graph has loops like that in Fig.3.

\begin{figure}
\centering
\includegraphics[height=3cm, width=10cm , bb=43 97 366 204]{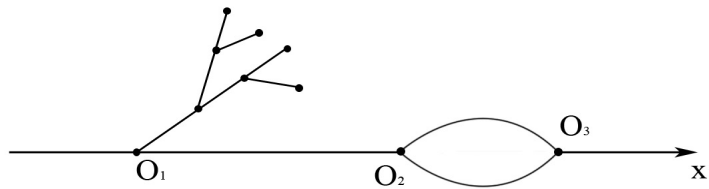}
\caption{A more general graph.}
\end{figure}

If the channel $D^\ve$ is not "uniformly narrow" but has points on
axis $x$ such that in the $\dt$-neighborhoods, $\dt=\dt(\ve)<\!<1$
of those points the channel has the "diameter" of order
$\mu(\ve)>\!>\ve$, the limiting process on the graph can have delays
or even traps. This will lead to different behavior of
$\dfrac{\sm^\ve_a}{a}$.

\end{document}